%% file: CPD.tex
\documentclass[conference,10pt]{IEEEtran}

\input{TeX/Preamble.tex}

	\shownotesfalse 
	\kfalse         
	\let\vec\bm
	\renewcommand{\norm}[1]{\|#1\|}

\begin{document}


\title{\sc A quadratically convergent proximal algorithm for nonnegative tensor
  decomposition\texorpdfstring{\thanks{\TheThanks}}{}}

\author{%
  Nico Vervliet\affil{1} \qquad Andreas Themelis\affil{1} \qquad Panagiotis
  Patrinos\affil{1} \qquad Lieven De Lathauwer\affil{1,2} \
  		\\[2.5mm]\small
		\affil{1}
			KU Leuven, Department of Electrical Engineering ESAT--STADIUS, Kasteelpark Arenberg 10, bus 2446, B-3001 Leuven, Belgium
		\\
		\affil{2}
			KU Leuven--Kulak, Group Science, Engineering and Technology, Etienne Sabbelaan 53, 8500 Kortrijk, Belgium
            \\
            {
			\{\href{mailto:Nico.Vervliet@kuleuven.be}{Nico.Vervliet},
			\href{mailto:Andreas.Themelis@kuleuven.be}{Andreas.Themelis},
			\href{mailto:Panos.Patrinos@kuleuven.be}{Panos.Patrinos},
			\href{mailto:Lieven.DeLathauwer@kuleuven.be}{Lieven.DeLathauwer\}@}\allowbreak
			\href{mailto:nico.vervliet@kuleuven.be,andreas.themelis@kuleuven.be,panos.patrinos@kuleuven.be,lieven.delathauwer@kuleuven.be}{kuleuven.be}%
		}%
      }

      \maketitle

	\begin{abstract}
      The decomposition of tensors into simple rank-1 terms is key in a variety
      of applications in signal processing, data analysis and machine
      learning. While this canonical polyadic decomposition (CPD) is unique
      under mild conditions, including prior knowledge such as nonnegativity can
      facilitate interpretation of the components. Inspired by the effectiveness
      and efficiency of Gauss--Newton (GN) for unconstrained CPD, we derive a
      proximal, semismooth GN type algorithm for nonnegative tensor
      factorization. If the algorithm converges to the global optimum, we show
      that $Q$-quadratic convergence can be obtained in the exact case. Global
      convergence is achieved via backtracking on the forward-backward envelope
      function. The $Q$-quadratic convergence is verified experimentally, and we
      illustrate that using the GN step significantly reduces number of
      (expensive) gradient computations compared to proximal gradient descent.
	\end{abstract}



	\section{Introduction}
		\input{TeX/Introduction.tex}

	\section{A semismooth Gauss-Newton method}\label{sec:GN}
		\input{TeX/GN.tex}

	\section{The forward-backward envelope}
		\input{TeX/FBE.tex}

	\section{A globally convergent algorithm}\label{sec:Algorithm}
		\input{TeX/Algorithm.tex}

   	\section{Complexity}\label{sec:complexity}
		\input{TeX/complexity.tex}

    \section{Experiments}\label{sec:experiments}
		\input{TeX/experiments.tex}

    \section{Conclusion and future work}\label{sec:conclusions}
		\input{TeX/conclusion.tex}


	\bibliographystyle{IEEEtranS}
	\bibliography{newbib.bib}

\end{document}

%% file: TeX/Preamble.tex
	\usepackage{algpseudocode}
	\usepackage{amsfonts}
	\usepackage{amsmath}
	\usepackage{amssymb}
	\usepackage{amsthm}
	\usepackage{array}
	\usepackage[USenglish]{babel}
	\usepackage{bm}
	\usepackage[nospace]{cite}
	\usepackage{comment}
	\usepackage{enumitem}
	\usepackage{etoolbox}
	\usepackage[T1]{fontenc}
	\usepackage{graphicx}
	\usepackage[utf8]{inputenc}
	\usepackage{mathtools}
	\usepackage{nicefrac}
	\usepackage{txfonts}
    \DeclareSymbolFont{symbols}{OMS}{cmsy}{m}{n}
	\usepackage{xargs}
	\usepackage[dvipsnames]{xcolor}
    \usepackage{extramath}
    \usepackage{pgf,pgfplots,pgfplotstable,tikz,advancedpictures}
	\usepackage[%
		breaklinks,
		plainpages    = false,
		colorlinks    = true,
		urlcolor      = BlueIMT,
		linkcolor     = MidnightBlue,
		citecolor     = Green,
		hypertexnames = false,
	]{hyperref}
	\usepackage[
		capitalize,
		nameinlink,
	]{cleveref}
	\usepackage{algorithm}

    \tikzexternaldisable
    \usetikzlibrary{positioning,calc,arrows,decorations}

	\newcommand{\affil}[1]{\texorpdfstring{\textsuperscript{#1}}{}}
	\newcommand{\TheThanks}{%
		This work was supported by the \emph{Research Foundation Flanders (FWO)}
        research projects G086518N, G086318N, and 12ZM220N; 
		\emph{Research Council KU Leuven} C1 projects No. C14/18/068 and C16/15/059;
		\emph{Fonds de la Recherche Scientifique---FNRS and the Fonds
          Wetenschappelijk Onderzoek---Vlaanderen} under EOS project No.
        30468160 (SeLMA). This research received funding from the Flemish
        Government under the “Onderzoeksprogramma Artificiële Intelligentie (AI)
        Vlaanderen” program.
	}

	\IEEEoverridecommandlockouts
	\definecolor{BlueIMT}{RGB}{0,42,72}
	\newcommand{\red}{\color{red}}                                      
	\newcommand{\blue}{\color{MidnightBlue}}                            
	\setlist[enumerate]{label={\it(\roman*)}}

	\newcommand\mathtight[1][0.25]{
		\thinmuskip  = #1\thinmuskip
		\medmuskip   = #1\medmuskip
		\thickmuskip = #1\thickmuskip
		\renewcommand\quad{\hspace{#1em}}%
		\renewcommand\qquad{\quad\quad}%
	}%
	\crefname{line}{step}{steps}
	\Crefname{line}{Step}{Steps}
	\renewcommand{\thealgorithm}{\Roman{algorithm}}
	\makeatletter
		\renewcommand{\theALG@line}{\thealgorithm.\oldstylenums{\arabic{ALG@line}}}
		\renewcommand{\alglinenumber}[1]{\theALG@line:}
	\makeatother

	\newif\ifk\ktrue                                                    
	\makeatletter
		\AtBeginDocument{%
			\ifk
				\gdef\@k+{k+1}%
				\gdef\k{\@ifnextchar+{\@k}{\@ifnextchar0{}{k}}}%
			\else
				\gdef\@k0{}%
				\gdef\k{\@ifnextchar0{\@k}{}}%
			\fi
		}
	\makeatother

	\newtheorem{thm}{Theorem}
		\newlist{enumthm}{enumerate}{1}
		\setlist[enumthm]{label={\it(\roman*)}, ref=\thethm{\it(\roman*)}}
		\crefalias{enumthmi}{thm}
		\AtBeginEnvironment{thm}{\let\enumerate\enumthm\let\endenumerate\endenumthm}
	\newtheorem{lem}[thm]{Lemma}
		\newlist{enumlem}{enumerate}{1}
		\setlist[enumlem]{label={\it(\roman*)}, ref=\thelem{\it(\roman*)}}
		\crefalias{enumlemi}{lem}
		\AtBeginEnvironment{lem}{\let\enumerate\enumlem\let\endenumerate\endenumlem}

	\DeclareMathOperator{\ball}{B}
	\DeclareMathOperator{\blockdiag}{blkdiag}
	\let\diag\relax
	\DeclareMathOperator{\diag}{diag}
	\DeclareMathOperator{\dist}{dist}
	\DeclareMathOperator{\id}{id}                                       
	\let\ker\relax
	\DeclareMathOperator{\ker}{null}                                    
	\DeclareMathOperator*{\minimize}{minimize}
	\DeclareMathOperator{\proj}{\Pi}
	\DeclareMathOperator{\range}{range}                                 
	\DeclareMathOperator{\stt}{subject\ to}
	\newcommand{\Fw}[1]{
		\ifstrempty{#1}{
			\id-\gamma JF^\top F
		}{
			#1-\gamma JF(#1)^\top F(#1)
		}
	}
	\newcommand{\FB}[1]{\proj_\C[\Fw{#1}]}                              

	\newcommand{\C}{\mathcal C}                                         
	\newcommand{\FBE}{\varphi_\gamma^{\text{\sc fb}}}                   
	\newcommand{\Id}{\ensuremath{\mat{I}}}                                         
	\newcommand{\N}{\varmathbb N}                                       
	\newcommand{\OO}{\mathbb O}                                         
	\newcommand{\R}{\varmathbb R}                                       
	\newcommand{\Res}{\mathcal R_\gamma}                                
	\renewcommand{\SS}{\mathbb S}                                       

	\newcommandx{\seq}[3][{2=k\in\N},{3={}}]{(#1)_{#2}^{#3}}            
	\newcommand{\innprod}[2]{\langle{}#1{},{}#2{}\rangle}               
	\newcommandx{\set}[2][{2={}}]{\{#1\ifstrempty{#2}{}{\mid#2}\}}      

\makeatletter
	\def\@command{}
	\def\@getscripts{%
		\let\@subarg\relax
		\let\@suparg\relax
		\@getscripts@testsub%
	}
	
	\def\@getscripts@testsub{\@ifnextchar_\@getscripts@sub\@getscripts@testsup}

	\def\@getscripts@testsup{\@ifnextchar^\@getscripts@sup\@getscripts@testsub@b}

	\def\@getscripts@sub_#1{%
		\def\@subarg{#1}%
		\ifx\@suparg\relax
			\expandafter\@getscripts@testsup@b
		\else
			\expandafter\@command
		\fi
	}

	\def\@getscripts@sup^#1{%
		\def\@suparg{#1}%
		\ifx\@subarg\relax
			\expandafter\@getscripts@testsub@b
		\else
			\expandafter\@command
		\fi
	}

	\def\@getscripts@testsub@b{%
		\ifx\@suparg\relax
			\def\@suparg{\@empty}
		\fi
		\@ifnextchar_\@getscripts@sub\@command%
	}

	\def\@getscripts@testsup@b{%
		\ifx\@subarg\relax
			\def\@subarg{\@empty}
		\fi
		\@ifnextchar^\@getscripts@sup\@command%
	}

	\let\oldoverbracket\overbracket
	\let\oldunderbracket\underbracket
	\let\oldoverbrace\overbrace
	\let\oldunderbrace\underbrace
	\newcommand\@overbracket[2][0.5pt]{%
		\oldoverbracket[#1]{#2}%
	}
	\newcommand\@underbracket[2][0.5pt]{%
		\oldunderbracket[#1]{#2}%
	}
	\def\underbracket{\@ifstar\@smashunderbracket\@underbracket}
	\def\overbracket{\@ifstar\@smashoverbracket\@overbracket}
	\def\underbrace{\@ifstar\@smashunderbrace\oldunderbrace}
	\def\overbrace{\@ifstar\@smashoverbrace\oldoverbrace}
	\newcommand{\@smashoverbracket}[2][0.5pt]{%
		\def\@optarg{#1}%
		\def\@arg{#2}%
		\def\@command{%
			\vphantom{\@arg}%
			\smash{\@overbracket[\@optarg]{\@arg}_{\@subarg}^{\@suparg}}%
		}%
		\@getscripts%
	}
	\newcommand{\@smashunderbracket}[2][0.5pt]{%
		\def\@optarg{#1}%
		\def\@arg{#2}%
		\def\@command{%
			\vphantom{\@arg}%
			\smash{\@underbracket[\@optarg]{\@arg}_{\@subarg}^{\@suparg}}%
		}%
		\@getscripts%
	}
	\newcommand{\@smashoverbrace}[1]{%
		\def\@arg{#1}%
		\def\@command{%
			\vphantom{\@arg}%
			\smash{\overbrace{\@arg}_{\@subarg}^{\@suparg}}%
		}%
		\@getscripts%
	}
	\newcommand{\@smashunderbrace}[1]{%
		\def\@arg{#1}%
		\def\@command{%
			\vphantom{\@arg}%
			\smash{\underbrace{\@arg}_{\@subarg}^{\@suparg}}%
		}%
		\@getscripts%
	}
\makeatother

	\newif\ifthisnotered
	\newif\ifshownotes\shownotestrue
	\newcommand{\note}[2][]{
		\ifshownotes
			\ifstrempty{#1}{\unskip}{}%
			\ifthisnotered
				\def\thisnotecolor{\red}%
				\global\thisnoteredfalse
			\else
				\def\thisnotecolor{\blue}%
				\global\thisnoteredtrue
			\fi
			{\thisnotecolor#1\rlap{\(^{\displaystyle\star}\)}}%
			\marginpar{\thisnotecolor%
				\smash{\rotatebox{90}{%
					\parbox{2cm}{\scriptsize\sl
						\hspace*{0pt}\llap{\(^{\displaystyle\star}\)}%
						#2%
					}%
				}}%
			}%
		\else
			\ifstrempty{#1}{\ignorespaces}{#1}%
		\fi
	}
	\AtEndPreamble{
		\ifshownotes
			\usepackage{mparhack}
		\fi
	}


%% file: TeX/Introduction.tex
The canonical polyadic decomposition (CPD) expresses an \(N\)th-order tensor
\(\ten{T}\) as a minimal number of \(R\) rank-1 terms, each of which is the
outer product, denoted by $\op$, of \(N\) nonzero vectors, with \(R\) the tensor
rank. Mathematically, we have
\begin{align}
  \label{eq:cpd}
	\ten{T}
{}={}
	\sum_{r=1}^R\vec{a}_r^{(1)}\op\cdots\op \vec{a}_r^{(N)} \eqqcolon \cp{\mat{A}^{(1)},\ldots,\mat{A}^{(N)}},
\end{align}
in which factor matrix $\mat{A}^{(n)}$ has $a_r^{(n)}$ as its columns.
The CPD is essentially unique under mild conditions, which is an attractive
property in many applications handling multiway data, e.g., in data analysis,
signal processing and machine learning
\cite{cichocki2015tensordecompositions,sidiropoulos2016tensordecompositions}.
To improve interpretability of the components, nonnegativity constraints are
often imposed on the factor vectors
\cite{cichocki2015tensordecompositions,sidiropoulos2016tensordecompositions},
i.e., \(\vec{a}_r^{(n)} \geq 0\), $n=1,\ldots,N$, $r=1,\ldots,R$, where the inequality is meant
elementwise.

The CPD can be cast as the following nonlinear least squares problem (NLS):
\[
	\minimize_{\vec{x}\in\R^d}
	f(\vec{x})\quad \text{with} \quad f(\vec{x})\coloneqq\tfrac12\|F(\vec{x})\|^2,
\]
where \(F(\vec{x})\) is a vector-valued polynomial (multilinear) function. The
Gauss--Newton method (GN) is a powerful tool to address this kind of problems,
as despite requiring only first-order information of \(F(\vec{x})\) it can exhibit up to
quadratic rates of convergence.  The idea behind GN is using the Gramian
\(JF(\vec{x})^\top JF(\vec{x})\) as a surrogate for \(\nabla^2f(\vec{x})\), which well
approximates the true Hessian around solutions \(\vec{x}^\star\) whenever
\(F(\vec{x}^\star)=\vec{0}\).  One iteration \(\vec{x}\mapsto \vec{x}^+\) of GN amounts to solving the
linear system
\[
	\bigl(JF(\vec{x})^\top JF(\vec{x})\bigr)(\vec{x}^+-\vec{x})
{}={}
	-JF(\vec{x})^\top F(\vec{x}),
  \]
  in which $JF(\vec{x})^\top F(\vec{x})$ is the gradient of $f(\vec{x})$.
  Thanks to the multilinear structure of the problem, the linear system can be
  solved efficiently using iterative methods
  \cite{sorber2013optimizationbasedalgorithms,vervliet2018numericaloptimization}.
  
  As is typical for higher-order methods, GN converges only if the starting
  point is already close enough to a solution, whence the need of a
  globalization strategy ensuring that the iterates eventually enter a basin of
  (fast) local convergence.  Thanks to the smoothness of the cost function
  \(f(\vec{x})\), many linesearch or trust region approaches can efficiently be
  employed for the purpose; see,
  e.g. \cite{sorber2013optimizationbasedalgorithms,vervliet2018numericaloptimization}.
  However, the nonsmoothness arising from the constraints makes the approach not
  applicable to \emph{nonnegative} CPD problems, namely
\begin{equation}\label{eq:P}
	\minimize_{\vec{x}\in\R^d}
	\tfrac12\|F(\vec{x})\|^2
\quad\stt{}\ 
	\vec{x}\geq0.
\end{equation}

In this paper we leverage on the globalization technique of
\cite{themelis2018forward,stella2017simple} to obtain a globally and
quadratically convergent algorithm for NCPD that directly addresses the
constrained formulation \eqref{eq:P}.  In fact, to further reduce the number of
singularities and nonoptimal stationary points, we impose an additional
nonconvex constraint that singles out ambiguities in the tensor decomposition
arising because of its equivalence up to scaling factors. We defer the details
to \cref{sec:GN}.

\subsection{Related work}
\label{sec:related-work}

A number of alternating least squares or block coordinate descent type methods
have been proposed to solve \cref{eq:P}. In these algorithms, one factor matrix
or one row or column is fixed at every iteration, after which a linear least
squares subproblem with nonnegativity constraints is solved
\cite{cichocki2009nonnegativematrix,cichocki2009fastlocal} by, e.g., using
multiplicative updates \cite{chi2012tensorssparsity}, active set methods
\cite{bro1998multiwayanalysis}, or the alternating direction method of
multipliers (ADMM) \cite{huang2016flexibleefficient}. To compute a nonnegative
CPD, other cost functions based on divergences can be used as well; see, e.g.,
\cite{chi2012tensorssparsity,hansen2015newtonbasedoptimization,cichocki2009nonnegativematrix}.

While these BCD methods are often easy to implement, their convergence is
slow. Therefore, a few algorithms based on GN or Levenberg--Marquardt (LM) have
been proposed. Nonnegativity constraints can then be enforced using logarithmic
penalty functions \cite{paatero1997weightednonnegative} or active set methods
\cite{sorber2013optimizationbasedalgorithms,vervliet2016tensorlab3,vervliet2018numericaloptimization,kelley1999iterativemethods}. By
change of variable, e.g., by replacing $x_i$ by $x_i^2$, \cref{eq:P} can be
converted to an unconstrained problem
\cite{royer2011computingpolyadic,sorber2015structureddatafusion}, which may lead
to a prohibitive increase of nonoptimal stationary points.
For nonnegative matrix factorization, a proximal LM type algorithm which solves
an optimization problem using ADMM in every iteration, has been proposed
\cite{huang2019globalsip}.

\subsection{Notation}
\label{sec:notation}

Scalars, vectors, matrices are denoted by lower case, e.g., \(a\), bold lower
case, e.g., \(\vec{a}\) and bold upper case, e.g., \(\mat{A}\), respectively.
Calligraphic letters are used for a tensor \(\ten{T}\), a constraint set $\C$, or
the uniform distribution $\ten{U}(a,b)$. Sets are indexed by superscripts within
parentheses, e.g., \(\mat{A}^{(n)}\), \(n=1,\ldots,N\). The Kronecker and
Khatri--Rao (column-wise Kronecker) products are denoted by \(\kron\) and
\(\kr\), respectively. The notation
$\blockdiag(\{\vec{x}^{(n)}_r\}_{r,n=1}^{R,N})$ is used for a block-diagonal
matrix with blocks
$\vec{x}_1^{(1)}, \vec{x}_2^{(1)}, \ldots, \vec{x}_R^{(1)}, \vec{x}_1^{(2)},
\ldots, \vec{x}_R^{(N)}$. The identity matrix is denoted by $\mat{I}$, the
column-wise concatenation of $\vec{a}$ and $\vec{b}$ by $[\vec{a};\vec{b}]$, and
the $\varepsilon$-ball around $\vec{x}$ by $\ball(\vec{x},\varepsilon)$.


%% file: TeX/GN.tex
We derive a GN method to compute the nonnegative CPD of an
$I_1\times I_2 \times \cdots \times I_N$ tensor $\ten{T}$. In order to prove
$Q$-quadratic convergence, the algorithm is developed for a slightly altered
problem in which the $(N-1)R$ degrees of freedom associated with the scaling
ambiguity are removed. We therefore require each vector to have unit norm and
explicitly isolate the magnitude of each term in the sum as a scalar, collected
in a vector \(\vec{\lambda}\in\R^R\) as
\[
	\ten{T}
{}={}
	\sum_{r=1}^R\lambda_r\cdot \vec{a}_r^{(1)}\op\cdots\op \vec{a}_r^{(N)}
\ \ \text{with }\ 
	\|\vec{a}_r^{(n)}\|=1, \forall n,r.
\]
The normalized version of problem \eqref{eq:P} thus becomes
\begin{equation}\label{eq:PN}
	\minimize_{\vec{a}\in\R^d,\vec{\lambda}\in\R^R}
	\tfrac12\|F(\vec{a},\vec{\lambda})\|^2
\ \stt{}
	\left\{
		\begin{array}{@{}l@{}}
			\vec{a},\vec{\lambda}\geq0,\\
			\|\vec{a}_r^{(n)}\|=1,
		\end{array}
	\right.
  \end{equation}
in which $F(\vec{a},\vec{\lambda}) = \sum_{r=1}^{R} \lambda_r\cdot
\vec{a}_r^{(1)}\op\cdots\op \vec{a}_r^{(N)} - \ten{T}$, $\vec{a} =
[\vec{a}^{(1)}_1;
\vec{a}^{(1)}_2;\ldots;\vec{a}^{(1)}_R,\vec{a}^{(2)}_1;\ldots;\vec{a}^{(N)}_{R}]$,
and $d=R\sum_{n=1}^{N}I_n$. 
The feasible set
\(
	\C
{}\coloneqq{}
	\set{(\vec{a},\vec{\lambda})\in\R^d\times\R^R}[\vec{a},\vec{\lambda}\geq0,\|\vec{a}_r^{(n)}\|=1]
\)
is nonconvex, but projecting onto it is a simple block-separable operation: one has
\(
	\proj_\C(\vec{a},\vec{\lambda})
{}={}
	(\tilde{\vec{a}},\tilde{\vec{\lambda}})
\)
where
\begin{equation}\label{eq:proj}
	\tilde{\vec{a}}_r^{(n)}
{}={}
	\proj_\SS\proj_{\OO_+}\vec{a}_r^{(n)}
{}={}
	\frac{[\vec{a}_r^{(n)}]_+}{\|[\vec{a}_r^{(j)}]_+\|}
\quad\text{and}\quad
	\tilde{\vec{\lambda}}
{}={}
	[\vec{\lambda}]_+.
\end{equation}
Here, \(\SS\) and \(\OO_+\) denote the unit sphere and the positive orthant of suitable size, respectively, and \([{}\cdot{}]_+=\max\set{0,{}\cdot{}}\) elementwise.
First-order necessary condition for optimality in this constrained minimization setting can be cast as the nonlinear equation \(\Res(\vec{x})=\vec{0}\), where \(\vec{x}\coloneqq(\vec{a},\vec{\lambda})\) is the optimization variable, and
\begin{equation}
	\Res(\vec{x})
{}\coloneqq{}
	\vec x-\proj_\C\bigl(\vec{x}-\gamma JF(\vec{x})^\top F(\vec{x})\bigr)
\end{equation}
is the projected-gradient residual mapping.
This map is everywhere piecewise smooth (up to a negligible set of points that we may disregard, as shown in the proof of \cref{thm:local}).
As such, its Clarke Jacobian \(J\Res\) furnishes a suitable first-order approximation.
The chain rule \cite[Prop. 7.1.11(a)]{facchinei2003finite} gives
\[\mathtight[0.25]
	J\Res(\vec{x})
{}={}
	\Id
	{}-{}
	J\proj_\C(\vec{w})\cdot\left[
		\Id
		{}-{}
		\gamma JF(\vec x){\hspace{-2pt}}^\top\hspace{-3pt} JF(\vec{x})
		{}-{}
		\underbracket*{
			\gamma\textstyle\sum_iF_i(\vec{x})\nabla^2F_i(\vec{x})
		}
	\right],
\]
where \(F_i(\vec{x})\) is the \(i\)-th element of vector \(F(\vec{x})\),
\begin{equation}\label{eq:Fw}
	\vec{w}
{}={}
	\vec{x}-\gamma JF(\vec{x})^\top F(\vec{x})
\end{equation}
is a gradient descent step at \(\vec{x}\), and \(J\proj_\C(\vec{w})\) is a (set of) \((d+R)\times (d+R)\) block-diagonal matrices.
In order to avoid Hessian evaluations, in the same spirit of (unconstrained) GN we replace \(J\Res\) with
\begin{equation}\label{eq:JR}
	\hat J\Res(\vec{x})
{}\coloneqq{}
	\Id
	{}-{}
	J\proj_\C(\vec{w})\cdot\left[
		\Id
		{}-{}
		\gamma JF(\vec{x})^\top JF(\vec{x})
	\right],
\end{equation}
which is \(O(\|\vec{x}-\vec{x}^\star\|)\)-close to \(J\Res(\vec{x})\) around a solution \(\vec{x}^\star\) of \eqref{eq:PN} provided that \(F(\vec{x}^\star)=\vec{0}\), as is apparent from the bracketed term in the expression of \(J\Res(\vec{x})\).
Since the feasible set \(\C\) is the product of small dimensional sets \(\SS_+\coloneqq\SS\cap\OO_+\) and \(\OO_+\), \(J\proj_\C\) is a structured set of block-diagonal matrices whose computation can be easily carried out using the chain rule
\(
	J\proj_{\SS_+}(\vec{w})
{}={}
	J\proj_\SS([\vec{w}]_+)
	J\proj_{\OO_+}(\vec{w})
\)
and the formulas
\begin{align}
  \label{eq:dproj}
  	J\proj_\SS([\vec{w}]_+)
{}={}
	\frac{\Id-\vec{z}\vec{z}^\top}{\|[\vec{w}]_+\|},
\text{with }
	\vec{z}
{}\coloneqq{}
	\proj_\SS([\vec{w}]_+)
{}={}
	\frac{[\vec{w}]_+}{\|[\vec{w}]_+\|},
\end{align}
and (see \cite[\S15.6.2d]{themelis2019acceleration})
\begin{align}
  \label{eq:jp+}
	J\proj_{\OO_+}(\vec{w})_{i,j}
	\begin{cases}
		=1 & \text{if }i=j \wedge w_i>0\\
		\in[0,1] & \text{if }i=j \wedge w_i=0\\
		=0 & \text{otherwise.}
	\end{cases}
\end{align}

\begin{thm}[Local quadratic convergence]\label{thm:local}%
	Let \(\vec{x}^\star\) be such that \(F(\vec{x}^\star)=\vec{0}\), and suppose that
	\note[all matrices in \(\hat J\Res(\vec{x}^\star)\) are nonsingular.]{Equivalent to strong local minimality?\\For \(\gamma\) small, yes...}
	Then, there exists \(\varepsilon>0\) such that the iterations
	\begin{equation}\label{eq:GN}
		\begin{cases}
			\vec{x}^0\in\ball(\vec{x}^\star,\varepsilon)\\
			\vec{x}^{k+1}=\vec{x}^k+\vec{d}^k,
		\end{cases}
	\quad\text{where}\quad
		\hat{\mat{H}}_k\vec{d}^k=-\Res(\vec{x}^k)
	\end{equation}
	with \(\hat{\mat{H}}_k\) being any element of \(\hat J\Res(\vec{x}^k)\), are \(Q\)-quadratically convergent to \(\vec{x}^\star\).
	\begin{proof}
		We start by remarking that the projection onto the (product of) sphere(s) is \(C^\infty\) wherever it is well defined.
		Since \(\vec{x}^\star\) is optimal, it follows from \cite[Thm. 3.4(iii)]{themelis2018forward} that \(\Res(\vec{x}^\star)=\set{\vec{0}}\), and that consequently the projection onto the spheres it entails, cf. \eqref{eq:proj}, is well defined and is thus \(C^\infty\) in a neighborhood.
		Combined with the strong semismoothness of the projection onto the positive orthant, see \cite[Prop. 7.4.7]{facchinei2003finite}, by invoking \cite[Prop. 7.4.4]{facchinei2003finite} we conclude that \(\Res\) is strongly semismooth around \(\vec{x}^\star\).

		Next, observe that \(\mat{H}_k\coloneqq
        \hat{\mat{H}}_k+\Delta(\vec{x}^k)\in J\Res(\vec{x}^k)\), for some
        \(\Delta(\vec{x})\in
        J\proj_\C(\vec{w})\sum_iF_i(\vec{x})\nabla^2F_i(\vec{x})\) (with
        \(\vec{w}\) as in \eqref{eq:Fw}) is a locally bounded quantity such that
        \(\Delta(\vec{x})\to \vec{0}\) as \(\vec{x}\to\vec{x}^\star\).
		Therefore, denoting \(\vec{d}^k\coloneqq \vec{x}^{k+1}-\vec{x}^k\),
		\[
			\|(\Res(\vec{x}^k)+\mat{H}_k)\vec{d}^k\|
		{}={}
			\|\Delta(\vec{x}^k)\vec{d}^k\|
		{}\leq{}
			\|\Delta(\vec{x}^k)\|\|\hat{\mat{H}}_k^{-1}\|\|\Res(\vec{x}^k)\|.
		\]
		We have that \(\sup_{\hat{\mat{H}}\in\hat J\Res(\vec{x})}\|\hat{\mat{H}}^{-1}\|\) is bounded by a same quantity \(c_\varepsilon\) for all \(\vec{x}\in\ball(\vec{x}^\star,\varepsilon)\) when \(\varepsilon\) is small enough, as it follows from \cite[Lem. 7.5.2]{facchinei2003finite}.
		Consequently, for \(\varepsilon\) small enough \cite[Thm. 7.5.5]{facchinei2003finite} guarantees that \(\vec{x}^k\to\vec{x}^\star\) \(Q\)-linearly.
		In turn, this implies that \(\Delta(\vec{x}^k)\to \vec{0}\), hence invoking again the same result the claimed \(Q\)-quadratic convergence is obtained.
	\end{proof}
\end{thm}

\Cref{thm:local} requires that all matrices in \(\hat J\Res(\vec{x}^\star)\)
(\cref{eq:JR}) are nonsingular. In \cref{thm:nonsingularHhat}, we show that this
is the case for the exact decomposition problem if the Gramian
$JF(\vec{x})^\top JF(\vec{x})$ has an $NR$-dimensional null space (which is
usually true for a unique CPD). This null space is derived in the next
lemma.

\begin{lem}[Kernel of Gramian]\label{lem:singularH}%
  The Gramian of the unconstrained problem $JF(\vec{x})^\top JF(\vec{x})$ has at
  least $NR$ zero eigenvalues, and a basis $\mat{K}$ for the subspace
  corresponding to these $NR$ zero eigenvalues is given by
  \begin{align}
    \mat{K} =
    \blockdiag \bigl(\{\diag (\vec{k}^{(n)})\kr\mat{A}^{(n)}\}_n,
    \diag (\vec{k}^{(N+1)})\kr \vec{\lambda}^{\top}\bigr),
  \end{align}
  for $\vec{k}^{(n)}\in\R^{R}$, $n=1,\ldots,N+1$, and \(\sum_{n=1}^{N+1}
  \vec{k}^{(n)} = \vec{0}\). 

  \begin{proof}
    It suffices to check that $JF(\vec{x}) \mat{K} = \mat{0}$ and that the
    dimension of $\mat{K}$ is $NR$. Using the expressions for $JF(\vec{x})$
    (see, e.g., \cite{vervliet2018numericaloptimization}) and multilinear
    identities, we have 
    \begin{align}
      JF(\vec{x})\mat{K} = \left(\kr_{n=1}^{N}
      \mat{A}^{(n)}\kr \vec{\lambda}^{\top}\right)\sum_{n=1}^{N+1}\vec{k}^{(n)}
      = \vec{0}.
    \end{align}
    As $\left(\kr_n \mat{A}^{(n)}\kr \vec{\lambda}^{\top}\right)$ usually has
    full column rank for an essentially unique decomposition defined by
    $\vec{x}$, we need \( \sum_{n=1}^{N+1} \vec{k}^{(n)} = \vec{0} \).  Since
    $\vec{k} = \shortmat{\vec{k}^{(1)}; \ldots; \vec{k}^{(N+1)}} \in
    \R^{(N+1)R}$ and the summation imposes $R$ linearly independent constraints,
    the columns of $\mat{K}$ span an $NR$-dimensional subspace.
  \end{proof}
\end{lem}

\begin{thm}\label{thm:nonsingularHhat}
  Let $\hat{\mat{H}}\in \hat{J}\Res(\vec{x}^{\star})$ and
  $F(\vec{x}^{\star}) = \vec{0}$. If the Gramian
  $JF(\vec{x}^{\star})^\top JF(\vec{x}^{\star})$ has $NR$ zero eigenvalues,
  $\hat{\mat{H}}$ is nonsingular.

    \begin{proof}
      In the global optimum $\vec{x}^{\star}$,
      $JF(\vec{x}^{\star})^{\T}F(\vec{x}^{\star})=\vec{0}$ and
      $\vec{x}^{\star}\in \C$, hence
      $J\Pi_{\C}(\vec{w}) = J\Pi_{\C}(\vec{x}^{\star})$. Let
      $\mat{P}\in J\Pi_{\C}(\vec{x}^{\star})$, and
      $\mat{G} = JF(\vec{x}^{\star})^\top JF(\vec{x}^{\star})$.  Before we prove
      that $\hat{\mat{H}}\in\hat{J}\Res(\vec{x}^{\star})$ has full rank, we show
      that $\range(\mat{P}\mat{G}) = \range(\mat{P})$ which is the case if
      $\range(\mat{P}) \cap \ker(\mat{G}) = \emptyset$. Let $\vec{a}^{(n)}_r$
      and $\lambda_r$ be the factor vectors and scaling factors corresponding to
      $\vec{x}^{\star}$. By assumption, $\mat{G}$ has $NR$ zero eigenvalues and
      $\ker(\mat{G})=\range(\mat{K})$; see \cref{lem:singularH}.  Any
      $\vec{b}\in\range(\mat{P})$ can be written as
      $\vec{b} = [\vec{b}^{(1)}_1, \ldots, \vec{b}^{(1)}_R, \vec{b}^{(2)}_1,
      \ldots, \vec{b}^{(N)}_R,b^{(N+1)}_1,\ldots,b^{(N+1)}_R]$ in which either
      $\vec{b}^{(n)}_r\perp \vec{a}_r^{(n)}$ or $\vec{b}^{(n)}_r = \vec{0}$; see
      \cref{eq:dproj}. If $\vec{b}\in \range(\mat{K})$, then the following
      should hold with \(\sum_{n=1}^{N+1} \vec{k}^{(n)} = \vec{0}\):
      \begin{align*}
        \vec{b}^{(n)}_r
        &= k^{(n)}_r \vec{a}_r^{(n)} \hspace{5.2mm}\Leftrightarrow
        & \hspace*{-1mm} k^{(n)}_r
        &= 0,
        & &\forall n,r,\\
        b^{(N+1)}_r
        &= k^{(N+1)}_r \lambda_r \hspace{4mm}\Leftrightarrow
        &\hspace*{-6mm} k^{(N+1)}_r
        &= b^{(N+1)}_r / \lambda_r,
        & &\forall r,
      \end{align*}
      which is false, hence $\vec{b}\notin \range(\mat{K})$ and
      $\range(\mat{P}\mat{G}) = \range(\mat{P})$.

      Let
      $\hat{\mat{H}} = \bigl(\mat{I} - \mat{P}\bigr) + \gamma \mat{P}\mat{G}$,
      and $l_0$ the number of active constraints for which
      $J\proj_{\OO_+}(\vec{x}^{\star})_{i,i} = 0$. As
      $\range(\mat{P}\mat{G})=\range(\mat{P})$, we can show that there exists an
      $NR+l_0$ dimensional subspace $U_1$ of $\mat{I} - \mat{P}$, and an $M-NR-l_0$
      dimensional subspace $U_2$ of $\mat{P}\mat{H}$, such that
      $\vec{u}_1^{\T}\vec{u}_2 = 0$, $\forall \vec{u}_{1}\in U_1$ and
      $\forall\vec{u}_2\in U_2$. Therefore,
      $\range(\bigl(\mat{I} - \mat{P}\bigr) + \gamma \mat{P}\mat{G}) = \R^{M}$ and
      $\hat{\mat{H}}$ has full rank.
    \end{proof}
\end{thm}


%% file: TeX/FBE.tex
\Cref{thm:local} highlights an appealing property that the constrained GN directions \eqref{eq:GN} enjoy close to the solutions of \eqref{eq:PN}.
Unfortunately, however, there is no practical way of initializing the iterations in such a way that the quadratic convergence is triggered.
In fact, not only is fast convergence not guaranteed without a proper initialization, but iterates may not converge at all and even diverge otherwise.
Because of the constraints, classical linesearch strategies cannot be adopted for \emph{nonnegative} CPDs.

Here, we overcome this limitation by integrating the fast GN directions \eqref{eq:GN} in the globalization strategy proposed in \cite{stella2017simple}, based on the \emph{forward-backward envelope} function \cite{patrinos2013proximal,themelis2018forward}
\begin{equation}\label{eq:FBE}
	\FBE(\vec x)
{}={}
	\tfrac12\|F(\vec x)\|^2
	{}-{}
	\innprod{JF(\vec x)^\top F(\vec x)}{\vec r}
	{}+{}
	\tfrac{1}{2\gamma}\|\vec r\|^2,
\end{equation}
where \(\gamma>0\) is a stepsize parameter,
\begin{equation}
	\vec z
{}\coloneqq{}
	\proj_\C\bigl(\vec x-\gamma JF(\vec x)^\top F(\vec x)\bigr),
\quad\text{and}\quad
	\vec r
{}\coloneqq{}
	\vec x-\vec z.
\end{equation}
The key properties of the FBE are summarized next.
Although an easy adaptation of that of \cite[Prop. 4.3 and Rem. 5.2]{themelis2018forward}, the proof is included for the sake of self containedness.

\begin{lem}[Basic properties of the FBE]\label{thm:FBE}%
	For every \(\gamma>0\), \(\FBE(\vec{x})\) is locally Lipschitz continuous and real-valued.
	Moreover, denoting \(f(\vec{x})\coloneqq\tfrac12\|F(\vec{x})\|^2\), the following hold: 
	\begin{enumerate}
	\item\label{thm:leq}
		\(\FBE(\vec x)\leq f(\vec x)\) for all \(\vec x\in\C\).
	\item\label{thm:geq}
		For all \(\vec x\in\R^d\), if \(\vec z\coloneqq\FB{\vec x}\) satisfies
		\begin{equation}\label{eq:UB}
			f(\vec z)\leq f(\vec x)+\innprod{\nabla f(\vec x)}{\vec z-\vec x}+\tfrac L2\|\vec x-\vec z\|^2
		\end{equation}
		for some \(L>0\), then
		\(
			F(\vec z)\leq\FBE(\vec x)-\tfrac{1-\gamma L}{2\gamma}\|\vec x-\vec z\|^2
		\).
		In particular, \(\FBE(\vec x)\geq f(\vec z)\geq0\) whenever \(\gamma\leq\nicefrac1L\).
	\item\label{thm:Lbounded}
		On every bounded set \(\Omega\subseteq\R^n\), there exists \(L_\Omega>0\) such that inequality \eqref{eq:UB} holds for every \(\vec x\in\Omega\) and \(L\geq L_\Omega\).
	\end{enumerate}
	\begin{proof}
		\renewcommand{\Fw}[1]{
			\ifstrempty{#1}{\id-\gamma\nabla f}{
				#1-\gamma\nabla f(#1)
			}
		}%
		Real valuedness is apparent from \cref{eq:FBE}.
		Moreover,
		\begin{align*}
			\FBE(\vec x)
		{}={} &
			\min_{\vec v\in\C}\set{f(\vec x)+\innprod{\nabla f(\vec x)}{\vec v-\vec x}+\tfrac12\|\vec x-\vec v\|^2}
		\\
		{}={} &
			f(\vec x)-\tfrac\gamma2\|\nabla f(\vec x)\|^2+\tfrac{1}{2\gamma}\dist\bigl(\vec x-\gamma\nabla f(\vec x),\,\C\bigr)^2
		\end{align*}
		hence \(\FBE(\vec x)\leq f(\vec x)\) whenever \(\vec x\in\C\) (by simply replacing \(\vec v=\vec x\) in the minimization).
		Local Lipschitz continuity owes to that of \(f\), \(\nabla f\) and \(\dist({}\cdot{},\C)\), see \cite[Ex. 9.6]{rockafellar2011variational}.
		Moreover, since the minimum above is obtained at \(\vec v=\vec z\), the second claim follows.
		Finally, since \(\proj_\C\) is locally bounded (cf. \cite[Ex. 5.23(a)]{rockafellar2011variational}) and so is \(\nabla f\), \(\FB{}\) maps the bounded set \(\Omega\) into a bounded set.
		Therefore, there exists a convex set \(\bar\Omega\) that contains all \(\vec x\) and \(\vec z=\FB{\vec x}\) with \(\vec x\in\Omega\).
		The claimed \(L_\Omega\) satisfying the last condition can thus be taken as the Lipschitz modulus of \(\nabla f\) over \(\bar\Omega\), see \cite[Prop. A.24]{bertsekas2016nonlinear}.
	\end{proof}
\end{lem}


%% file: TeX/Algorithm.tex
\begin{algorithm}
	\caption{PANOC for NCPD}
	\label{alg:PANOC}
	\input{TeX/PANOC.tex}
\end{algorithm}

\Cref{thm:FBE} contains all the key properties that lead to \cref{alg:PANOC}, which simply amounts to PANOC algorithm \cite{stella2017simple} specialized to this setting.
Having shown the efficacy of the fast GN directions \eqref{eq:GN}, the following result is a direct consequence of the more general ones in \cite{stella2017simple,themelis2018forward}.
We remark that the differentiability assumptions of \(\Res\) therein are only needed for showing the efficacy of quasi-Newton directions, whereas acceptance of unit stepsize only requires strong local minimality as shown in \cite[Thm. 5.23]{themelis2018proximal}%
\ifshownotes
	, \note[a condition entailed by the nonsingularity assumptions of \cref{thm:local}]{Elaborate}%
\fi
.

\begin{thm}[Convergence of \cref{alg:PANOC}]%
	Suppose that the sequence of points \(\vec z^{\k}\) remains bounded (as can be enforced by intersecting \(\C\) with any large box, cf. \cref{eq:tildeC}), then \Cref{alg:PANOC} terminates in finitely many iterations.
	Moreover, with tolerance \(\varepsilon=0\) the following hold:
	\begin{enumerate}
	\item
		\(\gamma\) is reduced only finitely many times and the sequence of points \(\vec z^{\k}\) converges to a stationary point for \eqref{eq:PN}.
	\item
		\note[If the conditions of \cref{thm:local} are satisfied at the limit point,]{Comment on nonsing. of Gramian}
		then eventually stepsize \(\tau_{\k}=1\) is always accepted and the sequence of points \(\vec z^{\k}\) converges \(Q\)-quadratically.
	\end{enumerate}
\end{thm}

Although the proof is already subsumed by previous work, in conclusion of this section we briefly outline the main details of the globalization strategy.
The algorithm revolves around the upper bound \eqref{eq:UB}; since the modulus \(L\) (initialized as \(L=\nicefrac\alpha\gamma\)) is not known a priori, it is adjusted adaptively throughout the iterations at \cref{step:test_gamma,step:test_gamma+} by halvening \(\gamma\) (hence doubling \(L=\nicefrac\alpha\gamma\) as a byproduct) until \eqref{eq:UB} is satisfied.
Since \(\vec z\) is the result of a projection on \(\C\), its \(\vec{\lambda}\)-component is nonnegative and its \(\vec{a}\)-component is bounded on unit spheres.
Consequently, boundedness of all the iterates may be artificially imposed by changing the feasible set \(\C\) into
\begin{equation}\label{eq:tildeC}
	\tilde\C
{}\coloneqq{}
	\set{(\vec{x},\vec{\lambda})\in\R^n\times\R^R}[\vec{a}\geq0,\|\vec{a}_r^{(n)}\|=1,0\leq\vec{\lambda}\leq
    M], 
\end{equation}
where \(M\) is a large constant.
Up to possibly resorting to this modification, as ensured by \cref{thm:Lbounded} the stepsize \(\gamma\) is halvened only a finite number of times and eventually remains constant.
Since \(\gamma=\nicefrac\alpha L<\nicefrac1L\), \cref{thm:geq} guarantees that
\[
	\tfrac12\|F(\vec z)\|^2
{}\leq{}
	\FBE(\vec x)
	{}-{}
	\tfrac{1-\alpha}{2\gamma}\|\vec r\|^2
{}<{}
	\FBE(\vec x)
	{}-{}
	\tfrac{1-\alpha}{2\gamma}\beta\|\vec r\|^2
\]
holds at every iteration; strict inequality holds because \(\beta<1\) and \(\vec r\neq\vec{0}\) (for otherwise the algorithm would have stopped at \cref{step:end}).
It then follows from the continuity of the FBE, \cref{thm:leq}, and the fact that \(\vec x^{\k+}\to\vec z\) as \(\tau_{\k}\searrow0\) (cf. \cref{step:x+}), that at every iteration \(\tau\) is halvened only a finite number of times at \cref{step:test_tau}.
In particular, the algorithm is well defined and, when \(\gamma\) becomes constant, produces a sequence satisfying
\[
	0
{}\leq{}
	\FBE(\vec x^+)
{}\leq{}
	\FBE(\vec x)
	{}-{}
	\tfrac{1-\alpha}{2\gamma}\beta\|\vec r\|^2.
\]
By telescoping the inequality, the vanishing of the \emph{residual} \(\vec r\) follows, hence the finite termination of the entire algorithm.


%% file: TeX/PANOC.tex
\begin{algorithmic}[1]
\def\mypar#1{\par\vspace{2pt}}
\small
\item[\sc\rlap{Require}\hphantom{Initialize}]
	\begin{tabular}[t]{@{}l@{}}
		Starting point \(\vec x^{\k0}\in\R^{d+R}\);~
		\(\alpha,\beta\in(0,1)\);~
		tolerance \(\varepsilon>0\);
	\\
		estimate of Lipschitz modulus \(L>0\)
	\end{tabular}
\item[\sc Initialize]
	\(\gamma=\nicefrac\alpha L\)%
	\ifk;~\(k=0\)\fi
\mypar{Preliminary tuning of \(\gamma\)}
\addtocounter{ALG@line}{-1}%
\State\label{step:test_gamma}%
	\(\vec z^{\k}=\FB{\vec x^{\k}}\)
	~and~
	\(\vec r^{\k}=\vec x^{\k}-\vec z^{\k}\)
\Statex{\bf if}~
	\(
		\tfrac12\|F(\vec z^{\k})\|^2
	{}>{}
		\FBE(\vec x^{\k})-\tfrac{1-\alpha}{2\gamma}\|\vec r^{\k}\|^2
	\)
	~{\bf then}
	\Statex\qquad
		\(\gamma\gets\nicefrac\gamma2\) and go back to \cref{step:test_gamma}
\mypar{Termination criterion}
\State\label{step:start}\label{step:end}%
	{\bf if}~
		\(\tfrac1\gamma\|\vec r^{\k}\|^2\leq\varepsilon\)
	~{\bf then}
\Statex\qquad
	{\bf return}~ \(\vec z^{\k}\)
\mypar{Computation of direction}
\State\label{step:direction}%
	Pick \(\hat{\mat{H}}_{\k}\in\hat J\Res(\vec x^{\k})\) and let \(\vec d^{\k}\) be such that
	{\ifshownotes
		\red
	\fi
	\(\hat{\mat{H}}_{\k}^\top\hat{\mat{H}}_{\k}\vec d^{\k}=-\hat{\mat{H}}_{\k}^\top\Res(\vec x^{\k})\)}%
\Statex
	Set stepsize \(\tau_{\k}=1\)
\mypar{Linesearch on \(\tau\) and tuning of \(\gamma\)}
\State\label{step:x+}%
	\(\vec x^{\k+} = (1-\tau_{\k})\vec z^{\k}+\tau_{\k}(\vec x^{\k}+\vec d^{\k})\)
\Statex
	\(\vec z^{\k+} = \FB{\vec x^{\k+}}\)
	~and~
	\(\vec r^{\k+} = \vec x^{\k+}-\vec z^{\k+}\)
\State{\bf if}~\label{step:test_gamma+}
	\(
		\tfrac12\|F(\vec z^{\k+})\|^2
		{}>{}
			\FBE(\vec x^{\k+})-\tfrac{1-\alpha}{2\gamma}\|\vec r^{\k+}\|^2
	\)
	~{\bf then}~
	\Statex\qquad
		\(\gamma\gets\nicefrac\gamma2\) and go back to \cref{step:test_gamma}
\Statex{\bf else if}~\label{step:test_tau}
	\(
		\FBE(\vec x^{\k+})
	{}>{}
		\FBE(\vec x^{\k})-\tfrac{1-\alpha}{2\gamma}{\beta}\|\vec r^{\k}\|^2
	\)
	\Statex\qquad
		\(\tau_{\k}\gets\nicefrac{\tau_{\k}}{2}\) and go back to \cref{step:x+}
\mypar{Proceed to next iteration}
\State\label{step:proceed}
	\ifk
		\(k\gets k+1\)
	\else
		\((\vec x^{\k}, \vec z^{\k}, \vec r^{\k})\gets(\vec x^{\k+}, \vec z^{\k+}, \vec r^{\k+})\)
	\fi
	and proceed to \cref{step:start}
\end{algorithmic}


%% file: TeX/complexity.tex
The computational complexity of \cref{alg:PANOC} is dominated by the same
operations as in the unconstrained case: computing the function value
$f(\vec{x})$ (\cref{step:test_gamma,step:test_gamma+}) and the gradient
$JF(\vec{x})^{\T}F(\vec{x})$ (\cref{step:test_gamma,step:x+}), and solving the
linear system \eqref{eq:GN} (\cref{step:direction}). If an iterative solver is
used to solve \cref{eq:GN}---which is common practice for medium to large scale
problems---the computational
complexity is dominated by the computation of the function evaluation and the
gradient, which require $O(R\prod_nI_n)$ and $O(NR\prod_nI_n)$ operations,
respectively, for an $N$th-order $I_1\times I_2 \times \cdots \times I_N$
tensor \cite{sorber2013optimizationbasedalgorithms}. Given a good initial
guess for the Lipschitz modulus $L$, \cref{step:test_gamma} is 
computed only a few times. Hence, the total complexity mainly depends on the
number of backtracking steps on $\tau$ at \cref{step:test_gamma+}.


%% file: TeX/experiments.tex
We validate the theoretical properties of \cref{alg:PANOC} by two
experiments. \Cref{alg:PANOC} is implemented in MATLAB 2019b with Tensorlab 3.0
\cite{vervliet2016tensorlab3}. Default values for the parameters
$\alpha=0.95$ and $\beta=0.5$ are used. If more than five backtracking steps on \(\tau\)
at \cref{step:test_gamma+} are needed, a proximal gradient
step is taken. The stopping tolerance is set to $\varepsilon=10^{-20}$ and the maximum number of iterations to 2000.
The Lipschitz modulus $L$ is estimated using
finite differences in a random direction. To prevent slower convergence due to
small $\gamma$, we heuristically set $\gamma\leftarrow \max\{\gamma,\eta\}$ in
\cref{step:proceed} with $\eta = \vec{g}^{\top}\mat{G}\vec{g}/\|\vec{g}\|^2$ in
which $\vec{g}$ and $\mat{G}$ are the gradient and Gramian for the unconstrained
problem, respectively. (The scaling factor $\eta$ is often used in trust region methods
to compute the Cauchy point.)

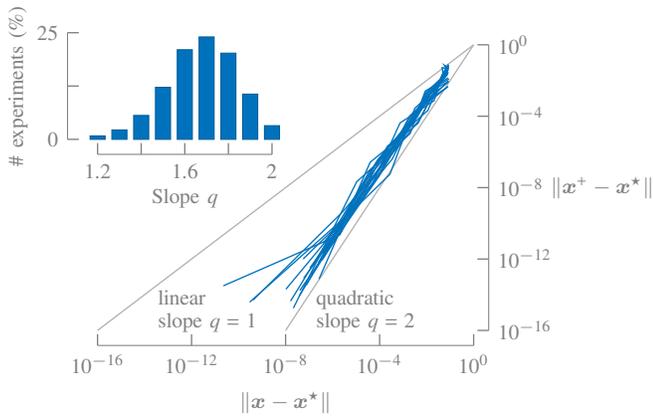
\begin{figure}[htp]
  \centering
  \input{TeX/experiment-quadratic.tikz}
  \caption{Up to $Q$-quadratic convergence can be achieved near the global optimum
    if an exact, unique solutions exists ($F(\vec{x}^{\star})=\vec{0}$). The histogram
    shows the slope for the penultimate iteration for 500 experiments. The
    convergence curves for ten randomly chosen experiments have a slope close to
    2. Results shown for a rank-5, $10\times10\times10$ tensor with a unique and
    nonnegative CPD, starting close to the global optimum.}
  \label{fig:experiment-quadratic}
\end{figure}

We show that $Q$-quadratic convergence can be achieved for exact nonnegative CPD. In
this experiment, 250 random $10\times10\times10$ tensors of rank-5 are
constructed using random factor matrices with entries drawn from
$\ten{U}(0,1)$. In each factor matrix, ten entries are set to zero at random to
ensure some constraints are active. For each random tensor, \cref{alg:PANOC} is
initialized by perturbing the exact solution such that approximately one digit
is correct. During the algorithm the distance to the exact solution
$\norm{\vec{x}^k - \vec{x}^{\star}}$ is tracked. This error should decrease as
$\|\vec{x}^{+} - \vec{x}^{\star}\| = O(\|\vec{x} - \vec{x}^{\star}\|^q)$
with $q=2$ to achieve $Q$-quadratic convergence, which is indeed the case, cf. \cref{fig:experiment-quadratic}.

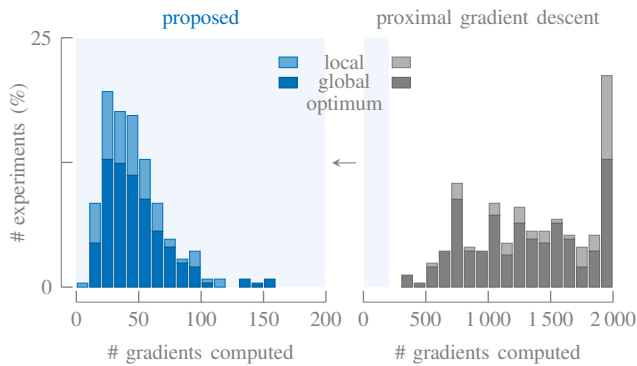
\begin{figure}[htp]
  \centering
  \input{TeX/experiment-global.tikz}
  \caption{By using (approximate) second-order information as in the proposed
    algorithm, fewer gradients are computed compared to proximal gradient
    descent. Histograms created using 250 experiments.}
  \label{fig:experiment-global}
\end{figure}

In the second experiment, we show that even for nonnegative tensor approximation
problems, the GN step can lead to a faster convergence. Similar to the
previous experiment, random $10\times10\times10$ tensors of rank-5 are
constructed using random factor matrices with entries drawn from
$\ten{U}(0,1)$. In each factor matrix, ten entries are replaced by small
negative entries drawn from $\ten{U}(-0.01,0)$. Hence, no exact nonnegative CPD
exists. Starting from a random initialization, the proposed method and standard
proximal gradient descent are run until convergence
($\norm{\vec{r}}^{2} < \gamma\varepsilon$). To eliminate excess iterations due to
nonoptimal stopping criteria, the number of gradient iterations is counted until
the algorithm converges to $1.01f(\vec{z}_{\text{final}})$, in which
$\vec{z}_{\text{final}}$ is the value returned by the algorithm. (This mainly
benefits proximal gradient descent.) As can be seen in
\cref{fig:experiment-global}, using the GN step clearly reduces the number of
gradient evaluations, which are the dominant cost. Note that both algorithms may
converge to local optima in which one or more of the rank-1 terms becomes zero.

%% file: TeX/experiment-quadratic.tikz
%
\definecolor{mycolor1}{rgb}{0.00000,0.44700,0.74100}%
\setlength{\figurewidth}{5cm}
\setlength{\figureheight}{3.8cm}

\begin{tikzpicture}[%
  join=round
  font=\footnotesize,
]

\begin{axis}[%
width=\figurewidth,
height=\figureheight,
at={(0\figurewidth,0\figureheight)},
scale only axis,
xmode=log,
xmin=1e-16,
xmax=1e0,
xminorticks=true,
xlabel={$\|\vec{x} -\vec{x}^\star\|$},
ymode=log,
ymin=1e-16,
ymax=1e-0,
ytick={1e0,1e-4,1e-8,1e-12,1e-16},
ylabel={$\|\vec{x}^{+} -\vec{x}^\star\|$},
minimal y right,
ylabel style={rotate=-90,text width=1.4cm,align=left,xshift=-2.5mm},
clip mode=individual,
every axis plot post/.style={line width=0.5pt,clip=true},
font=\footnotesize,
]
\addplot [color=minimal axis color,auxiliary line,forget plot]
  table[row sep=crcr]{%
1	1\\
1e-16	1e-16\\
};
\label{addplot:experiment-quadratic0}
\addplot [color=minimal axis color,auxiliary line,forget plot]
  table[row sep=crcr]{%
1	1\\
1e-08	1e-16\\
};
\label{addplot:experiment-quadratic1}
\addplot [color=mycolor1,solid,forget plot]
  table[row sep=crcr]{%
0.0845906837030552	0.052777104035632\\
0.052777104035632	0.0170097485817519\\
0.0170097485817519	0.00258539451416001\\
0.00258539451416001	0.000111260571062704\\
0.000111260571062704	1.36143448128433e-07\\
1.36143448128433e-07	3.75910887313603e-13\\
};
\label{addplot:experiment-quadratic2}
\addplot [color=mycolor1,solid,forget plot]
  table[row sep=crcr]{%
0.081722539260675	0.0799036070029271\\
0.0799036070029271	0.0426801454709061\\
0.0426801454709061	0.00875839940833269\\
0.00875839940833269	0.00257595517327378\\
0.00257595517327378	3.53619662011893e-05\\
3.53619662011893e-05	7.72772425601443e-08\\
7.72772425601443e-08	3.56665860100308e-14\\
};
\label{addplot:experiment-quadratic3}
\addplot [color=mycolor1,solid,forget plot]
  table[row sep=crcr]{%
0.0815664282593875	0.0141290977104075\\
0.0141290977104075	0.00305018193487985\\
0.00305018193487985	8.605675407228e-05\\
8.605675407228e-05	3.89477707202444e-07\\
3.89477707202444e-07	2.22893665355047e-12\\
};
\label{addplot:experiment-quadratic4}
\addplot [color=mycolor1,solid,forget plot]
  table[row sep=crcr]{%
0.08073069813252	0.00453518628763593\\
0.00453518628763593	0.000220454008098637\\
0.000220454008098637	3.52877040890634e-07\\
3.52877040890634e-07	8.57723967736464e-13\\
};
\label{addplot:experiment-quadratic5}
\addplot [color=mycolor1,solid,forget plot]
  table[row sep=crcr]{%
0.0849191256623212	0.0409645939162839\\
0.0409645939162839	0.0112782926287448\\
0.0112782926287448	0.000211394100847484\\
0.000211394100847484	2.03301473485364e-06\\
2.03301473485364e-06	2.2660301815693e-11\\
2.2660301815693e-11	3.16507011452027e-14\\
};
\label{addplot:experiment-quadratic6}
\addplot [color=mycolor1,solid,forget plot]
  table[row sep=crcr]{%
0.0837188975128957	0.00807831335098926\\
0.00807831335098926	0.000720439304124135\\
0.000720439304124135	2.58020538200376e-06\\
2.58020538200376e-06	4.18157619795369e-10\\
4.18157619795369e-10	5.10967818070407e-15\\
};
\label{addplot:experiment-quadratic7}
\addplot [color=mycolor1,solid,forget plot]
  table[row sep=crcr]{%
0.0801082194181155	0.0665950707377295\\
0.0665950707377295	0.049509360293243\\
0.049509360293243	0.0259705284342312\\
0.0259705284342312	0.00427796304765216\\
0.00427796304765216	3.95736373827387e-05\\
3.95736373827387e-05	2.64023974703982e-07\\
2.64023974703982e-07	7.55873756266731e-14\\
};
\label{addplot:experiment-quadratic8}
\addplot [color=mycolor1,solid,forget plot]
  table[row sep=crcr]{%
0.0828636555836828	0.0573897329666739\\
0.0573897329666739	0.01960835056913\\
0.01960835056913	0.00407229465426714\\
0.00407229465426714	9.3003239309972e-05\\
9.3003239309972e-05	2.6427565766861e-07\\
2.6427565766861e-07	5.855371628315e-13\\
};
\label{addplot:experiment-quadratic9}
\addplot [color=mycolor1,solid,forget plot]
  table[row sep=crcr]{%
0.0895377161378731	0.012281745399979\\
0.012281745399979	0.00266799991968463\\
0.00266799991968463	7.46918286936529e-05\\
7.46918286936529e-05	2.46366685743965e-07\\
2.46366685743965e-07	6.1507622092393e-13\\
};
\label{addplot:experiment-quadratic10}
\addplot [color=mycolor1,solid,forget plot]
  table[row sep=crcr]{%
0.0834776790651004	0.0648213433032728\\
0.0648213433032728	0.0174987319994641\\
0.0174987319994641	0.00277950361582697\\
0.00277950361582697	9.10098640871291e-05\\
9.10098640871291e-05	4.00463829598358e-08\\
4.00463829598358e-08	1.72489971389026e-14\\
};
\label{addplot:experiment-quadratic11}
\addplot [color=mycolor1,solid,forget plot]
  table[row sep=crcr]{%
0.0837560364262334	0.0516227140923719\\
0.0516227140923719	0.0165361372293252\\
0.0165361372293252	0.00205667483884403\\
0.00205667483884403	9.94184893784819e-05\\
9.94184893784819e-05	1.07283180242321e-07\\
1.07283180242321e-07	3.1107965557595e-13\\
};
\label{addplot:experiment-quadratic12}
\addplot [color=mycolor1,solid,forget plot]
  table[row sep=crcr]{%
0.0849229360186317	0.0097491985791878\\
0.0097491985791878	0.000567425386225773\\
0.000567425386225773	4.34273152737383e-06\\
4.34273152737383e-06	2.97339318975056e-10\\
2.97339318975056e-10	3.84698344682868e-15\\
};
\label{addplot:experiment-quadratic13}
\addplot [color=mycolor1,solid,forget plot]
  table[row sep=crcr]{%
0.0866296346435371	0.00905310397199269\\
0.00905310397199269	0.00174997278439662\\
0.00174997278439662	1.8875436791132e-05\\
1.8875436791132e-05	1.62313832244811e-08\\
1.62313832244811e-08	4.39365533832851e-15\\
};
\label{addplot:experiment-quadratic14}
\addplot [color=mycolor1,solid,forget plot]
  table[row sep=crcr]{%
0.0832995192202116	0.0753874580372361\\
0.0753874580372361	0.045223631471917\\
0.045223631471917	0.051860793172226\\
0.051860793172226	0.0223639646910095\\
0.0223639646910095	0.00189724601480994\\
0.00189724601480994	0.000103223827054258\\
0.000103223827054258	5.86570460669703e-08\\
5.86570460669703e-08	8.36190980477409e-14\\
};
\label{addplot:experiment-quadratic15}
\addplot [color=mycolor1,solid,forget plot]
  table[row sep=crcr]{%
0.0823090699354572	0.0754963525988797\\
0.0754963525988797	0.0565946318872777\\
0.0565946318872777	0.0718777262567005\\
0.0718777262567005	0.0465196119483908\\
0.0465196119483908	0.0157189485964115\\
0.0157189485964115	0.00236065781531124\\
0.00236065781531124	7.26674332305183e-05\\
7.26674332305183e-05	2.68673082677409e-07\\
2.68673082677409e-07	6.08027468823157e-13\\
};
\label{addplot:experiment-quadratic16}
\addplot [color=mycolor1,solid,forget plot]
  table[row sep=crcr]{%
0.0796012150142336	0.0227405516622508\\
0.0227405516622508	0.00312288236628037\\
0.00312288236628037	0.000274002769595319\\
0.000274002769595319	5.60529218228648e-08\\
5.60529218228648e-08	1.00705041351423e-12\\
};
\label{addplot:experiment-quadratic17}
\addplot [color=mycolor1,solid,forget plot]
  table[row sep=crcr]{%
0.0813986815674225	0.00460240626613484\\
0.00460240626613484	0.000338839311235232\\
0.000338839311235232	3.89338900487921e-07\\
3.89338900487921e-07	7.7577705504019e-12\\
};
\label{addplot:experiment-quadratic18}
\addplot [color=mycolor1,solid,forget plot]
  table[row sep=crcr]{%
0.0840549491742611	0.0290349861056447\\
0.0290349861056447	0.00442185352802675\\
0.00442185352802675	5.30151959299445e-05\\
5.30151959299445e-05	3.82452275633207e-08\\
3.82452275633207e-08	8.7242440658198e-15\\
};
\label{addplot:experiment-quadratic19}
\addplot [color=mycolor1,solid,forget plot]
  table[row sep=crcr]{%
0.0835192472696236	0.0465188530898284\\
0.0465188530898284	0.0206563059420127\\
0.0206563059420127	0.000738222505616275\\
0.000738222505616275	4.17605091846332e-05\\
4.17605091846332e-05	1.0290024937278e-08\\
1.0290024937278e-08	2.03422109678826e-14\\
};
\label{addplot:experiment-quadratic20}
\addplot [color=mycolor1,solid,forget plot]
  table[row sep=crcr]{%
0.0803002138863751	0.0296697423923985\\
0.0296697423923985	0.00878528476928988\\
0.00878528476928988	0.00198668314713757\\
0.00198668314713757	1.06142190510952e-05\\
1.06142190510952e-05	2.14269039893718e-08\\
2.14269039893718e-08	1.78812033839161e-15\\
};
\label{addplot:experiment-quadratic21}

\coordinate (hist) at (axis cs:1e-16,5e-6);

\node at (axis cs:8e-8,	 12e-16)	[anchor=west,text width=3cm] {quadratic\\slope $q=2$};
\node at (axis cs:15e-15,12e-16) [anchor=west,text width=3cm] {linear\\slope $q=1$};

\end{axis}

\pgfplotstableread{ 
  Center  Slope
  1.2000  0.0080 
  1.3000  0.0220 
  1.4000  0.0560 
  1.5000  0.1220 
  1.6000  0.2100 
  1.7000  0.2400 
  1.8000  0.2020 
  1.9000  0.1060 
  2.0000  0.0320 
}\slopetable

\begin{axis}[
  width=3.9cm,
  height=3cm,
  at=(hist),
  ybar,
  ymin=0,
  ymax=0.25,
  bar width={2mm},
  ytick={0, 0.125, 0.25},
  xtick={1.2,1.4,...,2},
  xticklabels={1.2,,1.6,,2},
  yticklabels={0,,25},
  xmin=1.2,
  xmax=2.0,
  minimal,
  name=slopehist,
  xlabel={Slope $q$},
  ylabel={\# experiments (\%)},
  ylabel style={inner sep=0pt,yshift=0mm},
  font=\footnotesize,
  minimal y axis offset = -4mm,
  xlabel style={yshift=1mm},
]
\addplot [fill=mycolor1,draw=mycolor1!80!black,opacity=1]
          table [y=Slope, x=Center] {\slopetable};   
\end{axis}

\end{tikzpicture}%

%% file: TeX/experiment-global.tikz
\definecolor{mycolor1}{rgb}{0.00000,0.44700,0.74100}%
\definecolor{mycolor2}{rgb}{0.85000,0.32500,0.09800}%
\colorlet{mycolor2}{gray}

\setlength{\figurewidth}{4.9cm}
\setlength{\figureheight}{4.9cm}

\begin{tikzpicture}[font=\footnotesize]
\pgfplotstableread{ 
Center  Local     Global
  5     0.0040         0  
 15     0.0400    0.0440  
 25     0.0680    0.1280  
 35     0.0520    0.1240  
 45     0.0600    0.1120  
 55     0.0400    0.0880   
 65     0.0280    0.0560   
 75     0.0080    0.0400   
 85     0.0040    0.0240   
 95     0.0160    0.0200   
105     0.0040    0.0040   
115     0.0080         0   
125          0         0   
135          0    0.0080   
145          0    0.0040   
155          0    0.0080   
165          0         0   
175          0         0   
185          0         0   
195          0         0   
}\panoctable

\pgfplotstableread{ 
Center  Local     Global
  50         0         0 
 150         0         0 
 250         0         0 
 350         0    0.0120 
 450         0    0.0040 
 550    0.0040    0.0200 
 650         0    0.0360 
 750    0.0160    0.0880 
 850    0.0040    0.0360 
 950         0    0.0360 
1050    0.0120    0.0720 
1150    0.0120    0.0320 
1250    0.0160    0.0640 
1350    0.0080    0.0480 
1450    0.0120    0.0440 
1550    0.0040    0.0640 
1650    0.0040    0.0480 
1750    0.0200    0.0200 
1850    0.0160    0.0360 
1950    0.0840    0.1280 
}\pgdtable

\coordinate (r1) at (0,0);
\coordinate (r2) at ($(r1)+(\figurewidth,\figureheight)-(16mm,16mm)$);

\fill [mycolor1!60!white,opacity=0.1] (r1) rectangle (r2);

\begin{axis}[
  width=\figurewidth,
  height=\figureheight,
  ybar stacked,   
  ymin=0,         
  ymax=0.25,
  bar width={1.4mm},
  ytick={0, 0.125, 0.25},
  xtick={0,50,...,200},
  yticklabels={0,,25},
  xmin=0,
  xmax=200,
  minimal,
  name=panoccpd,
  xlabel={\# gradients computed},
  ylabel={\# experiments (\%)},
  ylabel style={inner sep=0pt,yshift=-1mm},
]
\addplot [fill=mycolor1,draw=mycolor1!80!black,opacity=1]
table [y=Global, x=Center] {\panoctable};   
\addplot [fill=mycolor1!60!white,draw=mycolor1!100!white]
table [y=Local, x=Center] {\panoctable};
\end{axis}

\begin{axis}[
  at=(panoccpd.east),
  anchor=west,
  xshift=5mm,
  width=\figurewidth,
  height=\figureheight,
  ybar stacked,   
  ymin=0,         
  ymax=0.25,
  bar width={1.4mm},
  ytick={0, 0.125, 0.25},
  xtick={0,500,1000,...,2000},
  xmin=0,
  xmax=2000,
  minimal,
  axis y line=none,
  xlabel={\# gradients computed},
  name=pgd,
]
\addplot [fill=mycolor2,draw=mycolor2!80!black,opacity=1]
table [y=Global, x=Center] {\pgdtable};   
\addplot [fill=mycolor2!60!white,draw=mycolor2!100!white]
table [y=Local, x=Center] {\pgdtable};

\coordinate (r3) at (axis cs:0,0);
\coordinate (r4) at (axis cs:200,0.25);


\end{axis}

\fill [mycolor1!60!white,opacity=0.1] (r3) rectangle (r4);
\coordinate (c) at ($(r1)!0.5!(r2)$);
\draw [<-,>=stealth,shorten <=1mm, shorten >=1mm,minimal axis color] (c-|r2) --
(c-|r3);

\node at (panoccpd.north) [text=mycolor1,anchor=south] {proposed};
\node at (pgd.north)      [text=mycolor2,anchor=south] {proximal gradient descent};

\node at ($(pgd.north)!0.5!(panoccpd.north)$) [yshift=-9mm,text=minimal axis color] (opt) {optimum};
\node [above=-2mm of opt, text width=8mm,align=center, text=minimal axis color] (glob) {global};
\node [above=-2mm of glob, text width=8mm,align=center, text=minimal axis color] (loc) {local};

\node [left=1mm of loc,text width=1mm,inner
sep=2pt,fill=mycolor1!60!white,draw=mycolor1!100!white] {};
\node [left=1mm of glob,text width=1mm,inner
sep=2pt,fill=mycolor1,draw=mycolor1!80!black] {};

\node [right=1mm of loc,text width=1mm,inner
sep=2pt,fill=mycolor2!60!white,draw=mycolor2!100!white] {};
\node [right=1mm of glob,text width=1mm,inner
sep=2pt,fill=mycolor2,draw=mycolor2!80!black] {};

\end{tikzpicture}

%% file: TeX/conclusion.tex
By combining nonnegativity and unit-norm constraints, a proximal Gauss--Newton
type algorithm is derived. Global convergence is achieved by backtracking the GN
step to the proximal gradient descent (PGD) step based on the forward-backward
envelope function. While $Q$-quadratic convergence is only shown for the global
optima in the case of an exact, essentially unique decomposition, the GN
directions effectively reduce the computational cost compared to PGD.

In the current work we focused on theoretical properties; large-scale
implementations are part of future work.
